\documentclass[proc]{edpsmath}
\usepackage{url}
\usepackage{amssymb}
\usepackage{graphicx}
\usepackage{mathrsfs}
\usepackage{epstopdf}
\usepackage{psfrag}
\usepackage{framed}
\usepackage[colorlinks=true,bookmarksnumbered=true, pdfauthor={Jérémie Bettinelli}]{hyperref}

\theoremstyle{plain}
\newtheorem{thm}{Theorem}

\newtheorem{prop}[thm]{Proposition}

\newtheorem{corol}[thm]{Corollary}
\newtheorem{defi}{Definition}

\theoremstyle{definition}

\newcommand{\tol}{\xrightarrow[n \to \infty]{(d)}}
\newcommand{\tolk}{\xrightarrow[k \to \infty]{(d)}}
\newcommand{\Pb}{\mathbb{P}}
\newcommand{\m}{\mathfrak m}
\newcommand{\q}{\mathfrak q}
\newcommand{\lab}{\mathfrak l}
\newcommand{\tr}{\mathfrak t}
\newcommand{\de}{\mathrel{\mathop:}\hspace*{-.6pt}=}
\newcommand{\X}{\mathcal X}
\newcommand{\dish}{\delta_\mathcal{H}}
\newcommand{\dGH}{d_{\operatorname{GH}}}

\def\Mone{\mathcal{M}_1}
\def\muISE{\mu_{\mathrm{ISE}}}
\def\muISEsh{\mu_{\mathrm{ISE}}^{\mathrm{shift}}}

\begin{document}

\title{Random maps}
\author{C\'eline Abraham}\address{Universit\'e Paris Sud; \href{mailto:celine.abraham@math.u-psud.fr}{\nolinkurl{celine.abraham@math.u-psud.fr}}; \nolinkurl{http://www.math.u-psud.fr/~abraham}.}
\author{J\'er\'emie Bettinelli}\address{CNRS \& Institut \'Elie Cartan de Lorraine; \href{mailto:jeremie.bettinelli@normalesup.org}{\nolinkurl{jeremie.bettinelli@normalesup.org}}; \nolinkurl{www.normalesup.org/~bettinel}.}
\author{Gwendal Collet}\address{\'Ecole polytechnique; \href{mailto:gcollet@lix.polytechnique.fr}{\nolinkurl{gcollet@lix.polytechnique.fr}}; \nolinkurl{http://www.lix.polytechnique.fr/Labo/Gwendal.Collet}.}
\author{Igor Kortchemski}\address{DMA, \'Ecole Normale Sup\'erieure; \href{mailto:igor.kortchemski@normalesup.org}{\nolinkurl{igor.kortchemski@normalesup.org}}; \nolinkurl{http://www.normalesup.org/~kortchem}.}
%
%
\begin{abstract}
This is a quick survey on some recent works done in the field of random maps.
\end{abstract}
%
%
%
\maketitle
\section*{Introduction}

The study of random maps has generated a growing interest in the past decade. Since the pioneering work of Chassaing and Schaeffer~\cite{chassaing04rpl}, many contributions came nourishing the field, climaxing with the works of Le~Gall~\cite{legall11ubm} and Miermont~\cite{miermont11bms}, who established the convergence of rescaled uniform plane quadrangulations.
This proceeding aims at giving an overview of a few of the most notable results in the field. We refer the interested reader for example to~\cite{LGM12Buzios} or to Gr\'egory Miermont's Saint-Flour notes for further reading about random maps in general, as well as to the original works presented below for specific details.

\begin{figure}[ht]
	\centering\includegraphics[width=.75\linewidth]{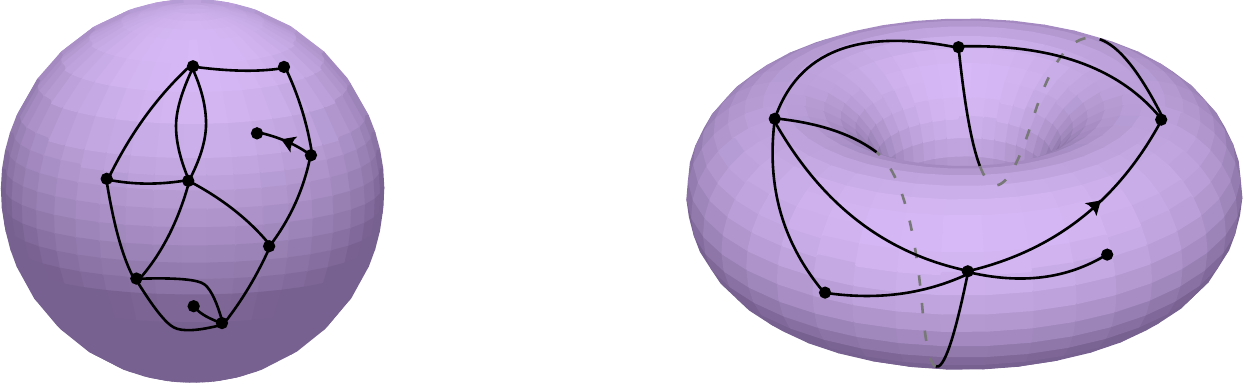}
	\caption{Two maps. The one on the left is a map of the sphere, the one on the right is a map of the torus.}
	\label{exmap}
\end{figure}

Let us first give a proper definition of the objects in question. Unfortunately, although quite natural and easy to comprehend, the concept of map is not very easy to define and the rigorous definition may seem a little bit complicated. We give it nonetheless for self-containment and refer to Figure~\ref{exmap} for an example. A \emph{map} is a cellular embedding of a finite graph (possibly with multiple edges and loops) into a compact connected orientable surface without boundary, considered up to orientation-preserving homeomorphisms. Cellular means that the faces of the map (the connected components of the complement of edges) are open $2$-cells, that is, homeomorphic to $2$-dimensional open disks. In particular, this implies that the graph is connected. In this paper, all the maps we consider are implicitely rooted, in the sense that they are given with a distinguished oriented edge called the \emph{root}.

There is also a more visual definition consisting in saying that a map is a gluing of polygons along their boundaries. The most studied maps are by far the maps of the sphere, also called \emph{plane maps} or \emph{planar maps}. In this case, the condition of cellular embedding is simply equivalent to the condition of connectedness of the graph. It is also natural to restrict one's attention to particular classes of maps. For example, one might consider \emph{triangulations}, that is, maps with faces of degree~$3$, where the \emph{degree} of a face is the number of oriented edges encountered when traveling along the boundary of the face. Beware that an edge may be visited twice in this contour, once in one direction, and once in the other one: such an edge will thus be counted twice. For example, on Figure~\ref{exmap}, the bottommost face of the map on the left has degree~$4$. We will see in the next section that \emph{quadrangulations} are also central objects (this is due to encoding combinatorial reasons to which we will come back in due course).
 
\bigskip

\textbf{Organization of the paper.} In the first section, we present a quick overview of the main results concerning limits of random maps. The subsequent sections summarize the original works presented during the random maps session of the \emph{Journ\'ees MAS 2014}.

\section{Limits of random maps}

\subsection{Scaling limits}\label{secsclim}

Maps come with a natural intrinsic metric given by the graph metric of the underlying graph. Let~$\m$ be a map and let $V(\m)$ denote its vertex-set. The \emph{graph metric} $d_\m$ is the metric on $V(\m)$ defined by
$$d_\m(u,v)\de\min\{k\ge 0\,:\, \text{ there exists a chain of $k$ edges linking $u$ to $v$}\},\qquad u,v\in V(\m).$$
Consequently, a map $\m$ may be seen as a finite metric space $(V(\m),d_\m)$ and it is natural to wonder what such a metric space looks like when the map is randomly distributed with a large number of faces. As the metric is growing with the number of faces, one has to rescale the space if one wants to observe a compact limit. This is the point of view of \emph{scaling limits}.

As we deal with isometry classes of metric spaces, the natural topology is the Gromov--Hausdorff topology introduced by Gromov~\cite{gromov99msr}. The Gromov--Hausdorff distance between two compact metric spaces $(\X,\delta)$ and $(\X',\delta')$ is defined by
$$\dGH\big((\X,\delta),(\X',\delta')\big) \de \inf \big\{ \dish\big(\varphi(\X),\varphi'(\X')\big)\!\big\},$$
where the infimum is taken over all isometric embeddings $\varphi : \X \to \X''$ and $\varphi':\X'\to \X''$ of~$\X$ and~$\X'$ into the same metric space $(\X'', \delta'')$, and $\dish$ stands for the usual Hausdorff distance between compact subsets of~$\X''$. This defines a metric on the set of isometry classes of compact metric spaces, making it a Polish space (see \cite{burago01cmg}).

The most natural setting is the following. We choose uniformly at random a map of ``size''~$n$ in some class, rescale the metric by the proper factor, and look at the limit in the sense of the Gromov--Hausdorff topology. The size considered is often the number of faces of the map. From this point of view, the most studied class is the class of plane quadrangulations. The pioneering work of Chassaing and Schaeffer~\cite{chassaing04rpl} revealed that the proper scaling factor in this case is~$n^{-1/4}$. The problem of scaling limit was then first addressed by Marckert and Mokkadem~\cite{marckert06limit}, who constructed a candidate limiting space called the \textit{Brownian map}, and showed the convergence toward it in another sense. Le~Gall~\cite{legall07tss} then showed the relative compactness of this sequence of metric spaces and that any of its accumulation points was almost surely of Hausdorff dimension~$4$. More precisely, he showed the following.

\begin{thm}[Le~Gall~\cite{legall07tss}]
Let $\q_n$ be uniformly distributed over the set of plane quadrangulations with~$n$ faces. Then, from any increasing sequence of integers, we may extract a subsequence $(n_k)_{k\ge 0}$ such that there exists a random metric space $(\q_\infty,d_\infty)$ satisfying
$$\Big( V(\q_{n_k}), n_k^{-1/4}\, d_{\q_{n_k}} \Big) \tolk \big(\q_\infty,d_\infty\big)$$
in the sense of the Gromov--Hausdorff topology.

Moreover, regardless of the choice of the sequence of integers, the Hausdorff dimension of the limiting space $(\q_\infty,d_\infty)$ is almost surely~$4$.
\end{thm}

 It is only recently that the solution of the problem was completed independently by Miermont~\cite{miermont11bms} and Le~Gall~\cite{legall11ubm}, who showed that the scaling limit is indeed the Brownian map. This last step, however, is not mandatory in order to identify the topology of the limit: Le~Gall and Paulin~\cite{legall08slb}, and later Miermont~\cite{miermont08sphericity}, showed that any possible limit is homeomorphic to the $2$-dimensional sphere. Summing up, we obtain the following theorems.

\begin{thm}[Le~Gall~\cite{legall11ubm}, Miermont~\cite{miermont11bms}]
Let $\q_n$ be uniformly distributed over the set of plane quadrangulations with~$n$ faces. There exists a random metric space $(M,D)$ called the \emph{Brownian map} such that
$$\Big( V(\q_{n}), \Big(\frac{9}{8n}\Big)^{1/4}\, d_{\q_{n}} \Big) \tol (M,D)$$
in the sense of the Gromov--Hausdorff topology.
\end{thm}

Note that the factor $(9/8)^{1/4}$ does not seem to be of importance in the previous statement but the proper definition of the Brownian map is with this factor. It comes from historical reasons and from an alternative description in terms of random processes.

\begin{thm}[Le~Gall~\cite{legall07tss}, Le~Gall--Paulin~\cite{legall08slb}, Miermont~\cite{miermont08sphericity}]
Almost surely, the Brownian map $(M,D)$ is homeomorphic to the $2$-dimensional sphere and has Hausdorff dimension~$4$.
\end{thm}

\begin{figure}[ht]
	\centering\href{https://sketchfab.com/models/aeabacc4b0ad49059f6b5a8e97ebd60c/embed?transparent=1}{\includegraphics[width=.95\linewidth]{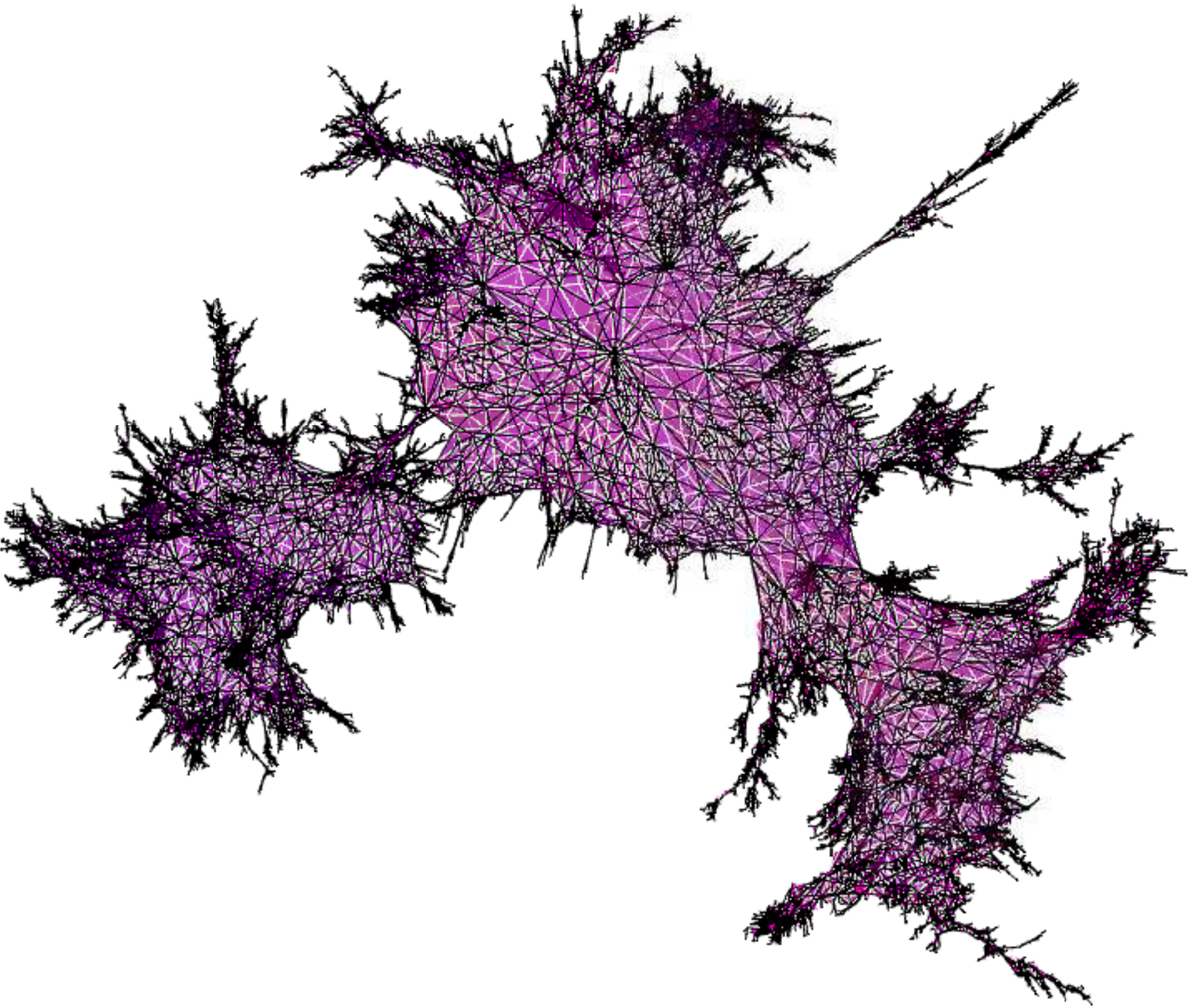}}
	\caption{Simulation of a uniformly sampled plane quadrangulation with $50\,000$ faces. See \url{http://www.normalesup.org/~bettinel/simul.html} for more computer simulations of random maps.}
	\label{sphere}
\end{figure}

This line of reasoning lead the way to several extensions. The first kind of extension is to consider other classes of plane maps. Actually, Le~Gall already considered in~\cite{legall07tss} the classes of $\kappa$-angulations, for even $\kappa\ge 4$. In~\cite{legall11ubm}, he considered the classes of $\kappa$-angulations for $\kappa=3$ and for even $\kappa\ge 4$ as well as the case of Boltzmann distributions on bipartite\footnote{A map is called \emph{bipartite} if its vertex-set may be partitioned into two disjoint subsets such that no edges link two vertices from the same subset.} plane maps, conditioned on their number of vertices. Another extension is due to Addario-Berry and Albenque~\cite{addarioalbenque2013simple}; in this work, they consider simple triangulations and simple quadrangulations, that is, triangulations and quadrangulations without loops and multiple edges. Bettinelli, Jacob and Miermont~\cite{BeJaMi14} later added the case of maps conditioned on their number of edges, and Abraham~\cite{abr13} considered the case of bipartite maps conditioned on their number of edges. In all these cases, the limiting space is always the same Brownian map (up to a multiplicative constant): we say that the Brownian map is \emph{universal} and we expect it to arise as the scaling limit of a lot more of natural classes of maps. A peculiar extension is due to Le~Gall and Miermont~\cite{legall09scaling} who consider maps with large faces, forcing the limit to fall out of this universality class: they obtain so-called \emph{stable maps}, which are related to stable processes.

Another kind of extension is to consider quadrangulations on a fixed surface that is no longer the sphere. The case of orientable surfaces with a boundary was the focus of~\cite{bettinelli14gbs}, in which the existence of subsequential limits was established. The complete convergence (without the need of extracting subsequences) is the subject of the work~\cite{bm14} in preparation:

\begin{thm}[Bettinelli--Miermont~\cite{bm14}]
Let $g\ge 0$ be an integer, $\sigma=(\sigma^1,\dots,\sigma^p)$ be a $p$-uple of positive real numbers and $\sigma_n=(\sigma_n^1,\dots,\sigma_n^p) \in \mathbb N^p$ be such that $\sigma_n^i/{\sqrt{2n}} \to \sigma^i$, for $1\le i \le p$. Let $\q_n$ be uniformly distributed over the set of bipartite genus~$g$-maps with~$n$ faces of degree~$4$ and~$p$ distinguished faces of respective degrees~$\sigma_n^1$, \ldots, $\sigma_n^p$. There exists a random metric space $(M^{(g)}_{\sigma},D)$ such that
$$\Big( V(\q_{n}),  \Big(\frac{9}{8n}\Big)^{1/4}\, d_{\q_{n}} \Big) \tol (M^{(g)}_{\sigma},D)$$
where the convergence holds in distribution for the Gromov--Hausdorff topology.

Moreover, almost surely, the limiting space $M^{(g)}_\sigma$ is homeomorphic to the surface of genus~$g$ with~$p$ boundary components, has Hausdorff dimension~$4$, and every of the~$p$ connected components of its boundary has Hausdorff dimension~$2$.
\end{thm}

\subsection{Geodesics toward a uniformly chosen point}

The metric of the Brownian map $(M,D)$, and a fortiori of more general Brownian surfaces $(M^{(g)}_{\sigma},D)$ is quite hard to grasp. Due to their construction, one can however study the geodesics from all the points to a given point uniformly distributed at random. In the Brownian map, Le~Gall~\cite{legall08glp} did this study and showed among other properties that, almost surely, the following holds.
\begin{itemize}
	\item There exists only one geodesic between two uniformly distributed points.
	\item Given a uniformly distributed point~$\rho$, there exist at most $3$ distinct geodesics from a point toward~$\rho$.
\end{itemize}

Moreover, he showed the following so-called confluence property of the geodesics.
\begin{prop}[Le~Gall~\cite{legall08glp}]
Given a uniformly distributed point~$\rho$, almost surely, for every $\varepsilon>0$, there exists $\eta\in\,(0,\varepsilon)$ such that all the geodesics from~$\rho$ to points outside of the ball of radius~$\varepsilon$ centered at~$\rho$ share a common initial part of length~$\eta$.
\end{prop}

The previous properties were later extended to every Brownian surface $(M^{(g)}_{\sigma},D)$ by Bettinelli~\cite{bettinelli14gbs}, by using a completely different approach. As, in the general case, the topology is richer than in the case of the sphere, some new interesting characterizations may be made in terms of geodesics. To give only one example, it is shown that, given a uniformly distributed point~$\rho$, there is only a finite number of points linked to~$\rho$ by three distinct geodesics such that the concatenation of any two such geodesics form a loop that is not homotopic to~$0$. This number is random and bounded by an explicit function of the genus and number of boundary components. Its distribution may also be computed. For instance, in the case of a surface of genus~$g\ge 1$ without boundary, the situation is simpler: this number is deterministic and equal to $4g-2$.

\subsection{Encoding by simpler objects: Schaeffer-like bijections}\label{secschaeffer}

The starting point of the previous studies is a powerful bijective encoding of the maps in the studied class by simpler objects. In the case of plane quadrangulations, the bijection in question is the so-called Cori--Vauquelin--Schaeffer bijection~\cite{cori81planar,schaeffer98cac,chassaing04rpl} and the simpler objects are so-called \emph{well-labeled trees}, that is, trees whose vertices carry integer labels such that the root vertex has label~$0$ and the variation of labels along the edges are in $\{-1,0,1\}$. In the other cases, variants of this bijection are used \cite{bouttier04pml,chapuy07brm,poulscha06,ambjornbudd} and the encoding objects usually have a more intricate combinatorial structure. In this section, we present the simplest of these bijections, that is, the original Cori--Vauquelin--Schaeffer bijection. 

\begin{figure}[ht]	
	\centering\includegraphics[width=.95\linewidth]{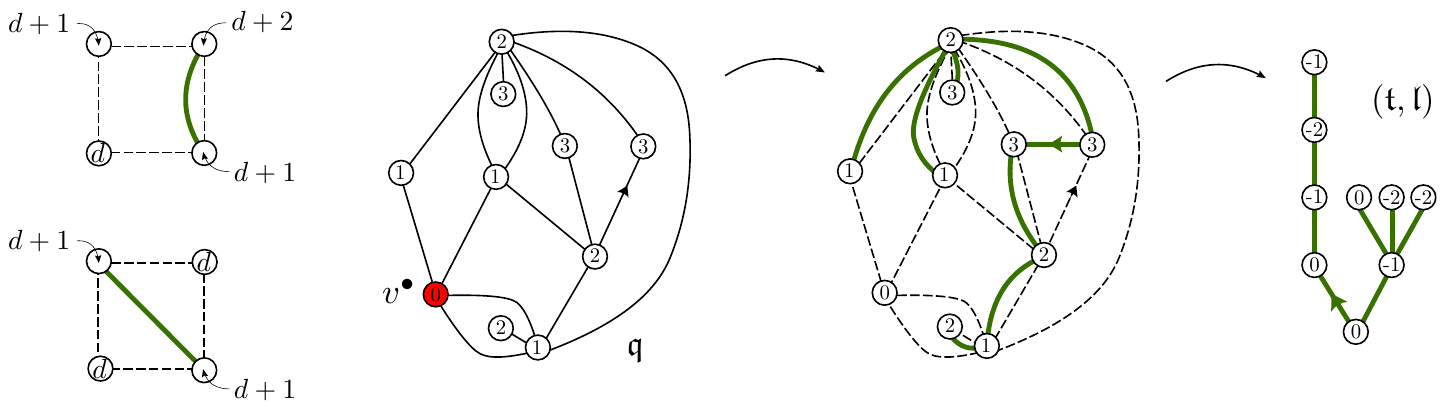}
	\caption{\textbf{Left.} The two types of faces and the extra arc to draw. \textbf{Right.} Construction of the well-labeled tree $(\tr,\lab)$ from the pointed quadrangulation $(\q,v^\bullet)$.}
	\label{cvsqt}
\end{figure}

\bigskip\noindent\textbf{From pointed plane quadrangulations to well-labeled trees.}
We start from a pointed quadrangulation $(\q,v^\bullet)$ with~$n$ faces. First, we assign to every vertex~$v$ a label $\hat\lab(v) \de d_\q(v,v^\bullet)$. Due to combinatorial reasons, the labels of neighboring vertices differ exactly by~$1$. As a result, when we travel around a face, the subsequent labels are either~$d$, $d+1$, $d+2$, $d+1$, or~$d$, $d+1$, $d$, $d+1$, for some integer~$d$. We draw inside each face an arc, as depicted on the left part of Figure~\ref{cvsqt}. The embedded graph with vertex-set $V(\q)\setminus\{v^\bullet\}$ and edge-set the set of newly added arcs is a tree~$\tr$, rooted according to some convention which we will not need to discuss in this work. Finally, the vertices of~$\tr$ are labeled by translating the function~$\hat\lab$ in such a way that the root vertex~$\rho$ of~$\tr$ gets label~$0$, that is, we set $\lab\de\hat\lab-\hat\lab(\rho)$. See Figure~\ref{cvsqt}.

\begin{figure}[ht]
	\centering\includegraphics[width=.95\linewidth]{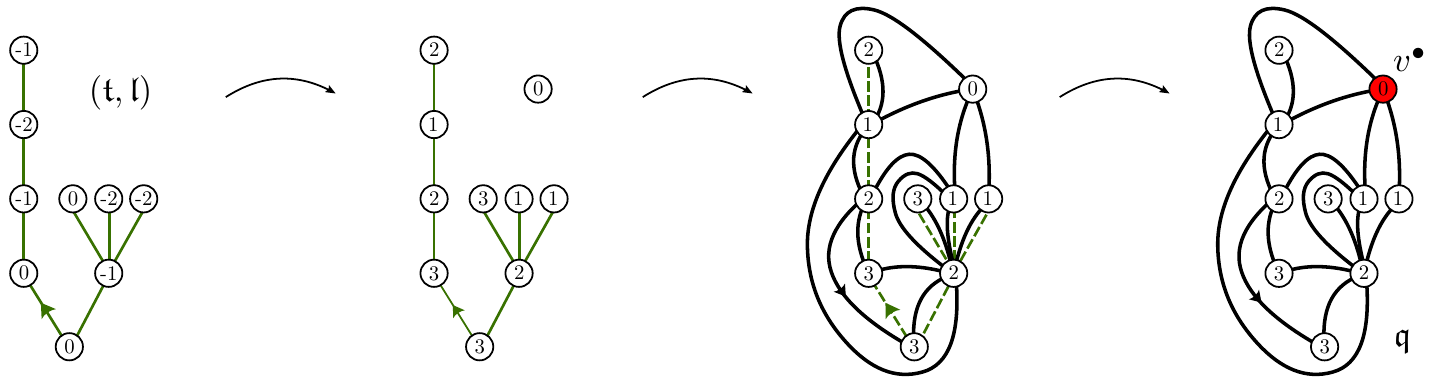}
	\caption{Reverse construction. We first translated the labels in such a way that the minimum label equals~$0$ in order to recover the labels from the previous section. This step does not influence the construction.}
	\label{cvstq}
\end{figure}

\bigskip\noindent\textbf{From well-labeled trees to pointed quadrangulations.}
Let $(\tr,\lab)$ be a well-labeled tree with~$n$ edges and let~$c_0$, $c_1$, \ldots, $c_{2n-1}$ denote the sequence of corners of~$\tr$, arranged in clockwise order around the tree, starting from the root corner. We extend the list to $\mathbb Z_+\cup\{\infty\}$ by periodicity, setting $c_{2kn+i}=c_{i}$ for every $k\in\mathbb N$, $i\in \{0,\ldots,2n-1\}$, and by adding one corner~$c_\infty$ incident to an extra vertex~$v^\bullet$ lying inside the unique face of~$\tr$. We also assign the label $\lab(v^\bullet)\de \min_{V(\tr)} \lab -1$. We then define the successor function $s:\mathbb Z_+\to\mathbb Z_+\cup\{\infty\}$ by
\begin{equation}\label{succ}
s(i)\de	\inf\{j>i:\lab(c_j)=\lab(c_i)-1\},\qquad i\in\mathbb Z_+,
\end{equation}
where, for a corner~$c$, we wrote~$\lab(c)$ the label of the incident vertex and, as usual, $\inf\{\varnothing\}\de\infty$. We then set $s(c_i)\de c_{s(i)}$. The construction consists in linking every corner~$c$ with its successor~$s(c)$ by an arc, in a non-crossing fashion. The embedded graph with vertex-set $V(\tr)\cup\{v^\bullet\}$ and edge-set the set of added arcs is then a quadrangulation, rooted according to some convention. See Figure~\ref{cvstq}.

The previous constructions are respectively $2$-to-$1$ and $1$-to-$2$, and are reverse one from another. The factor~$2$ comes from the rooting convention: the two quadrangulations obtained by reversing the root correspond to the same well-labeled tree.

\subsection{The Brownian snake}\label{secbs}

The first step in the comprehension of scaling limits of random maps is often the comprehension of the scaling limits of the encoding objects. In the case of a uniform plane quadrangulation with~$n$ faces, the encoding object is a uniform well-labeled tree with~$n$ edges and the scaling limit of such an object is a function of the \emph{Brownian snake}, introduced by Le~Gall~\cite{legall99sbp}. 

First, let~$\mathbf e=(\mathbf e_t)_{0\le t\le 1}$ be a normalized Brownian excursion (it can for instance be obtained by normalizing the excursion including time~$1$ of a reflected standard Brownian motion). Then, the so-called \emph{Brownian snake's head} driven by~$\mathbf e$ may be defined as the process $(\mathbf e_t,Z_t)_{0\le t \le 1}$, where, conditionally given~$\mathbf e$, the process~$Z$ is a centered Gaussian process with covariance function
\begin{equation*}
\operatorname{Cov}\!\big(Z_s, Z_t\big) = \inf_{s\wedge t\le u\le s\vee t} \mathbf e_u.
\end{equation*}

\begin{figure}[ht]
	\centering\includegraphics[width=.95\linewidth]{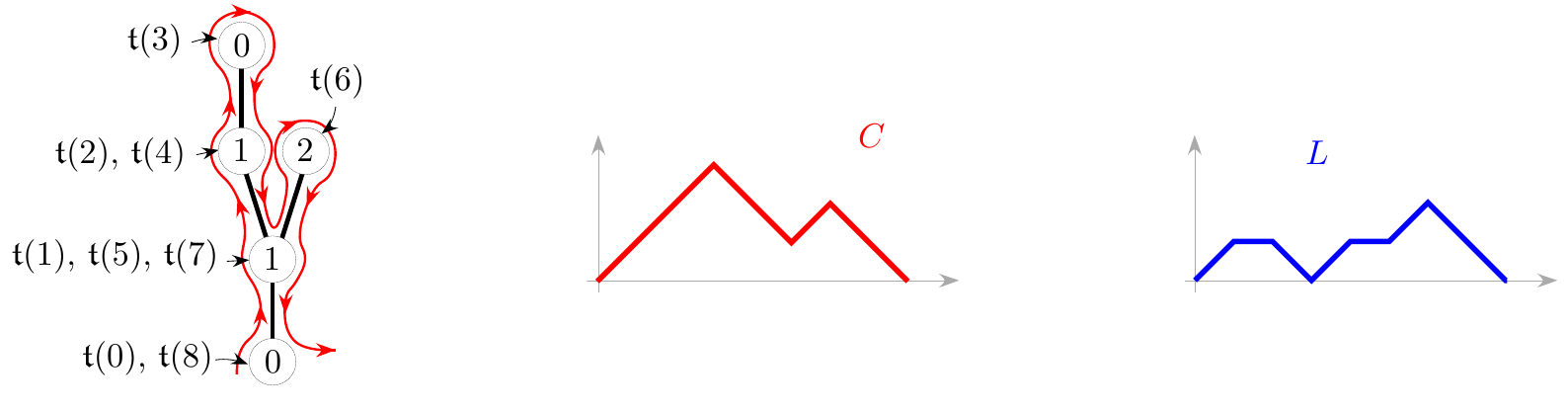}
	\caption{The contour and label function associated with a well-labeled tree having~$4$ edges.}
	\label{contour}
\end{figure}

Let us encode a well-labeled tree $(\tr,\lab)$ with~$n$ edges by its so-called contour and label functions~$C$ and~$L$ defined as follows. First, let $\tr(0)$, $\tr(1)$, \dots, $\tr(2n)$ be the vertices of~$\tr$ read in clockwise order around the tree, starting at the root corner. The \emph{contour function} $C:[0,2n]\to \mathbb R_+$ and the \emph{label function} $L:[0,2n]\to \mathbb R$ are defined by
$$C(i) \de d_{\tr}\big(\tr(0),\tr(i)\big)\qquad\text{ and }\qquad L(i) \de \lab(\tr(i)),\qquad 0 \le i \le 2n,$$
and linearly interpolated between integer values (see Figure~\ref{contour}).

\begin{thm}[Chassaing--Schaeffer \cite{chassaing04rpl}]\label{cvcontour}
Let~$C_n$ and~$L_n$ be the contour and label functions associated with a uniform well-labeled tree having~$n$ edges. Then
$$\Big(\Big(\frac1{2n}\Big)^{1/2} \,C_n(2n\,t), \, \Big(\frac{9}{8n}\Big)^{1/4}\, L_n(2n\,t)\Big)_{0 \le t \le 1}  
				\tol (\mathbf e_t,Z_t)_{0\le t \le 1}$$
for the uniform topology on $\mathcal C( [0,1],\mathbb R^2)$.
\end{thm}

The idea is roughly the following. The convergence of the first coordinate is a conditioned version of Donsker's theorem, which is due to Kaigh~\cite{kaigh76ipr}. Then, conditionally given the tree~$\tr_n$, the labels are given by i.i.d.\ uniform random variables on $\{-1,0,1\}$, which have variance equal to~$2/3$. As a result, the label variation between times~$i$ and~$j>i$ is given by a random walk of length~$d_{\tr_n}(\tr(i),\tr(j))$, which is of order $\sqrt{2n}\, (\mathbf e_{i/2n} + \mathbf e_{j/2n} -2\min_{i\le k\le j} \mathbf e_{k/2n})$ as a consequence of the convergence of the first coordinate. Donsker's theorem yields that this random walk approximates in the scale $ (9/8n)^{1/4}$ a Brownian motion of length $\mathbf e_{i/2n} + \mathbf e_{j/2n} -2\min_{i\le k\le j} \mathbf e_{k/2n}$.

\bigskip
We end this section by defining the one-dimensional ISE (random) probability measure $\muISE$ by
$$\muISE(h) \de \int_0^1 dt \, h(Z_t)$$
for every nonnegative measurable function~$h$.

\subsection{Local limits}\label{secll}

In the previous sections, we rescaled the metric in order to obtain a compact metric space at the limit. Another well-studied point of view consists in looking at \emph{local limits}, giving noncompact discrete objects.

In this setting, both triangulations and quadrangulations are often studied. Let~$\tr_n$ be a triangulation chosen uniformly at random among all triangulations with~$n$ vertices. Angel and Schramm~\cite{angel03uipt} have introduced an infinite random plane triangulation $T_{\infty}$, called the \emph{Uniform Infinite plane Triangulation (UIPT)}, which is obtained as the local limit of $\tr_{n}$ as $ n \to \infty$. More precisely, $T_{\infty}$ is characterized by the fact that, for every $r \geq 0$, we have the following convergence in distribution
$$B_{r}(\tr_{n}) \quad \xrightarrow[n\to\infty]{(d)} \quad B_{r}(T_{\infty}),$$
where $B_r(\m)$ denotes the map formed by the edges and vertices of the map $\m$ lying at graph distance smaller than or equal to~$r$ from the origin of the root.  This infinite random triangulation and its quadrangulation analog (the UIPQ, see \cite{chassaing06uipq,krikun05uipq}) have attracted a lot of attention, see \cite{BenjaCur13,curien2013glimpse,GurelNach13} and the references therein. We will also come back to this object in Section~\ref{seckor}, where percolation on it is studied.
  
\bigskip
Note that the scaling factor appearing in the results of Section~\ref{secsclim} is always $n^{-1/4}$ and, in this section, we did not rescale the metric. One might wonder what happens if one rescales by a different factor, say $n^{-\alpha}$ for some $\alpha\in(0,1/4)$. The resulting object is called the \emph{Brownian plane} and was introduced by Curien and Le~Gall~\cite{curien2012brownian}. This object may also be obtained by scaling the metric of the UIPQ by a factor~$\lambda$ and by letting $\lambda\to 0$.

\bigskip
\begin{framed}
\centering\itshape
The subsequent sections summarize the original works presented during the random maps session of the \emph{Journ\'ees MAS 2014}. 
\end{framed}

\section{Rescaled bipartite plane maps converge to the Brownian map [by C. Abraham]}

The work~\cite{abr13} is an addition to a series of papers that focus on the convergence of large random plane maps viewed as metric spaces to the continuous random metric space known as the Brownian map. The main goal is to provide another example of these limit theorems, in the case of bipartite plane maps with a fixed number of edges. 

Recall that a map is \emph{bipartite} if its vertices can be colored with two colors in such a way that any two vertices that have the same color are not connected by an edge. In the case of plane maps, this is equivalent to the property that all faces of the map have an even degree. Let $\mathbf{M}^b_n$ stand for the set of all rooted bipartite plane maps with~$n$ edges. 
\begin{thm}[Abraham \cite{abr13}]
\label{cvgcarte} 
For every $n\geq 1$,
let $\m_n$ be uniformly distributed over $\mathbf{M}^b_n$. Then, the following convergence holds in distribution for the Gromov--Hausdorff topology,
$$\big(V(\m_n) , (2n)^{-1/4}\, d_{\m_n}\big)  \tol (M, D)$$
where $(M, D)$ is the Brownian map introduced during Section~\ref{secsclim}.
\end{thm}

As the number of bipartite plane maps with any fixed number of vertices is always infinite, it would make no sense 
to condition on the number of vertices, and for this reason we consider conditioning
on the number of edges, which results in certain additional technical difficulties. 

The proof of Theorem~\ref{cvgcarte} relies on a generalization of the combinatorial bijection presented in Section~\ref{secschaeffer}, which is due to Bouttier, Di~Francesco and Guitter~\cite{bouttier04pml}. It specializes into a bijection between rooted and pointed bipartite maps and certain two-type labeled trees. The random labeled tree associated with a uniform rooted and pointed bipartite map with $n$ edges via this bijection is identified as a labeled two-type Galton--Watson tree with explicit offspring distributions, conditioned to have a fixed progeny $n$. This tree is denoted by $\big(\mathcal{T}_n, (\ell_n(u))_{u \in \mathcal{T}_n^0}\big)$, where $\mathcal{T}_n^0\subseteq V(\mathcal T_n)$ is the set of labeled vertices of $\mathcal{T}_n$. 
In order to prove the convergence to the Brownian map, an important technical step is to derive asymptotics for the contour and label functions associated with the conditioned tree $\big(\mathcal{T}_n, (\ell_n(u))_{u \in \mathcal{T}_n^0}\big)$.

\begin{thm}[Abraham \cite{abr13}]
Let $\big(\mathcal{T}, (\ell(u))_{u \in \mathcal{T}^0}\big)$  be a two-type Galton--Watson tree as previously, but with no condition on the total progeny. Let $N=|V(\mathcal{T})|-1$ be the random size of $\mathcal{T}$. 
The conditional distribution of $\big(n^{-1/2}\, C^{\mathcal{T}^0}(n\,t),\allowbreak n^{-1/4}\, L^{\mathcal{T}^0}(n\,t) \big)_{0 \leq t \leq 1}$ knowing that $N=n$ converges as $n$ tends to infinity to the law of $\big( \frac{4\sqrt{2}}{9}\, \mathbf{e}_t, 2^{1/4}\, Z_t \big)_{0 \leq t \leq 1}$, where $(\mathbf{e},Z)$ was introduced during Section~\ref{secbs}. 
\end{thm}

The fact that the tree is conditioned on the total number of vertices instead of the number of vertices of one type creates an additional difficulty, which we handle through an absolute continuity argument. A useful technical ingredient is a seemingly new definition of a ``modified'' Lukasiewicz path associated with a two-type tree.

The fact that the encoding bijection uses \emph{pointed} maps also introduces another new technicality. Indeed, note that a uniform pointed plane quadrangulation with~$n$ faces is a uniform plane quadrangulation with~$n$ faces, which is uniformly pointed. This is because every plane quadrangulation with~$n$ faces has $n+2$ vertices (this is an easy consequence of the Euler characteristic formula). In the case of uniform bipartite plane maps, this fact no longer holds and we first obtain the convergence of uniform pointed plane bipartite maps. We then derive Theorem~\ref{cvgcarte} at the cost of a technical ``depointing'' lemma.

\section{Percolation on random triangulations \& looptrees [by I. Kortchemski]}\label{seckor}

The work \cite{CKpercolooptrees} studies site percolation on Angel \& Schramm's UIPT (see Section~\ref{secll}). Its goal is to compute several critical and near-critical exponents, and describe the scaling limit of the boundary of large percolation clusters in all regimes (subcritical, critical and supercritical). One of the contributions of~\cite{CKpercolooptrees} is to prove in particular that the scaling limit of the boundary of large critical percolation clusters is the random stable looptree of index $3/2$, which was introduced in \cite{CK13}. 

The probabilistic theory of random plane maps and its physics counterpart, the Liouville 2D quantum gravity, is a very active field of research. The goal is to understand universal large-scale properties of random plane graphs or maps.
One possible way to get information about the geometry of these random lattices is to understand the behavior of (critical) statistical mechanics models on them. In~\cite{CKpercolooptrees}, Curien and Kortchemski focus on one of the simplest of such models: site percolation on the UIPT, which we now introduce.

\bigskip\noindent\textbf{Percolation on the UIPT.} Given  the UIPT, we consider a site percolation by coloring its vertices independently white with probability $a\in (0,1)$ and black with probability $1-a$. This model has already been studied by Angel \cite{angel03gpu}, who proved that the critical threshold  is almost surely equal to $$a_{c} =   {1}/{2}.$$ His approach was based on a clever Markovian exploration of the UIPT called the \emph{peeling process}. See also \cite{ACpercopeel,curien2013glimpse,MenardNolin14} for further studies of  percolation on random maps using the peeling process.

\medskip

\begin{figure}[!ht]
 \begin{center}
 \includegraphics[width=0.8 \linewidth]{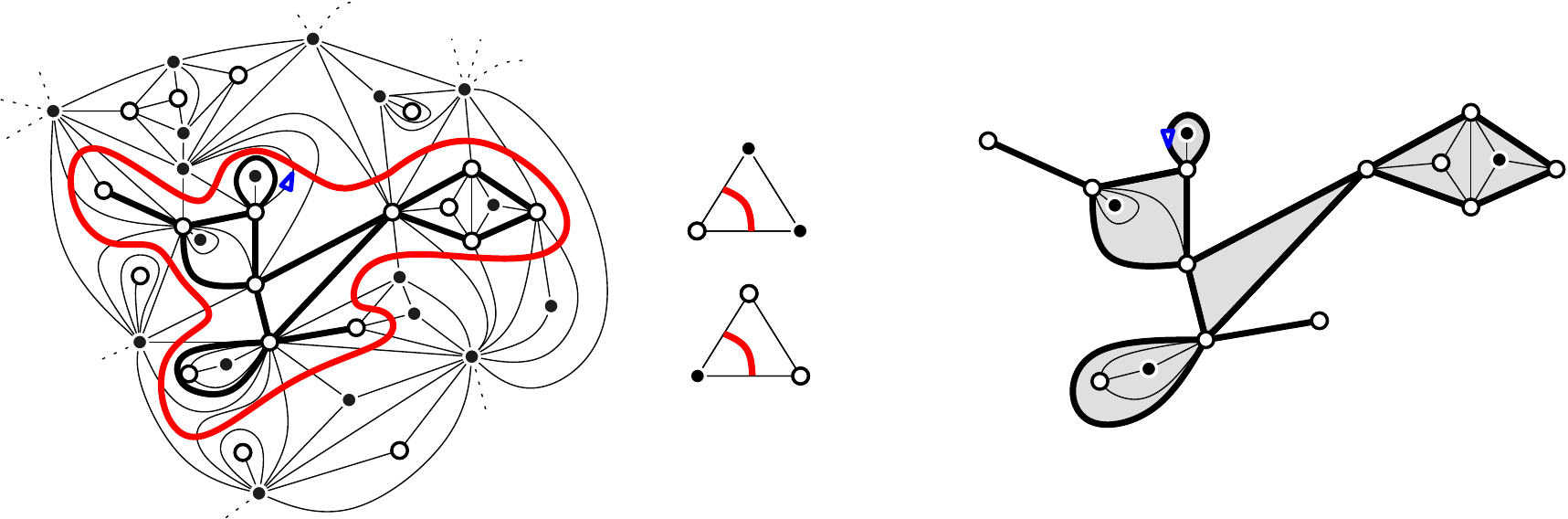}
 \caption{\label{fig:clusterloop}On the left, a part of a site-percolated triangulation with the interface going through the root edge and the hull of the cluster of the origin. The interface is drawn using the rules displayed in the middle triangles. On the right, the boundary of the hull is in bold black line segments and has perimeter $16$.\label{percoex-site}} 
 \end{center}
 \end{figure}

In~\cite{CKpercolooptrees}, the authors are interested in the geometry of \emph{the boundary of percolation clusters} and use a different approach. We condition the root edge of the UIPT on being of the form $\circ \to \bullet$, which will allow us to define the percolation interface going through the root edge. The white cluster of the origin is by definition the set of all the white vertices and edges between them that can be reached from the origin of the root edge using white vertices only.  We denote by $ \mathcal{H}^{\circ}_{a}$ its hull, which is obtained by filling-in the holes of the white cluster except the one containing the target of the root edge called the exterior component (see Figure~\ref{percoex-site}). Finally,  we denote by  $ \partial  \mathcal{H}_{a}^\circ$ the boundary of the hull, which is the graph formed by the edges and vertices of  $\mathcal{H}^{\circ}_{a}$ adjacent to the exterior (see Figure~\ref{percoex-site}), and let $ \#\partial  \mathcal{H}_{a}^\circ$ be its perimeter, or length, that is the number of oriented edges of $\partial  \mathcal{H}_{a}^\circ$ belonging to the exterior. Note that  $  \partial \mathcal{H}^{\circ}_{a}$ is formed of discrete cycles attached by some pinch-points.  It follows from the work of Angel \cite{angel03gpu} that, for every value of $a \in (0,1)$, the boundary $\partial \mathcal{H}_{a}^\circ$ is always \emph{finite} (if an infinite interface separating a black cluster from a white cluster existed, this would imply the existence of both  infinite black and white clusters, which is intuitively not posible).

One of the contributions of~\cite{CKpercolooptrees} is to find the precise asymptotic behavior for the probability of having a large perimeter in the critical case.

 \begin{thm}[Curien--Kortchemski \cite{CKpercolooptrees}] \textbf{\textup{Critical exponent for the perimeter.}} \label{thm:expo}

For $a= a_{c}=1/2$ we have 
 $$ \Pb\left( \# \partial \mathcal{H}^{ \circ}_{\frac{1}{2}} = n\right) \quad  \underset{n \to \infty}{\sim} \quad \frac{ {3}}{ {2} \cdot |\Gamma(-{2}/{3})|^3} \cdot n^{-4/3},$$
 where $ \Gamma$ is Euler's Gamma function.
 \end{thm}

It is interesting to mention that the exponent $4/3$ for the perimeter of the boundary of critical clusters also appears when dealing with the half-plane model of the UIPT: using the peeling process, it is shown in \cite{ACpercopeel} that $ P( \# \partial \mathcal{H}_{a_{c}}^\circ > n) \asymp n^{-1/3}$, where  $a_{n} \asymp b_{n}$ means that the sequence $a_{n}/b_{n}$ is bounded from below and above by certain constants.

One of the techniques used to establish Theorem \ref{thm:expo} is a surgery operation inspired from Borot, Bouttier \& Guitter \cite{BBG12,BBG12more}. This allows the authors to derive a tree representation of the $2$-connected components of $  \partial \mathcal{H}_{a}^{\circ}$, which is proved to be closely related to the law of a certain two-type Galton--Watson tree. The study of this two-type random tree is reduced to the study of a standard one-type Galton--Watson tree by using a recent bijection due to Janson \& Stef\'ansson \cite{JS12},  which enables the use of the vast literature on random trees and branching processes to make exact computations.

This method also allows to fully understand the probabilistic structure of the hull of the white cluster and to identify the scaling limits (for the Gromov--Hausdorff topology) in any regime (subcritical, critical and supercritical)  of $\partial \mathcal{H}_{a}^{\circ}$, seen as a compact metric space, when its perimeter tends to infinity. In particular, the authors establish that the scaling limit of $\partial \mathcal{H}_{a_{c}}^{\circ}$ conditioned to be large, appropriately rescaled, is the stable looptree of parameter $ 3/2$ introduced in \cite{CK13}, whose definition we now recall.

\begin{figure}[!ht]
 \begin{center}
 \includegraphics[width=0.65  \linewidth]{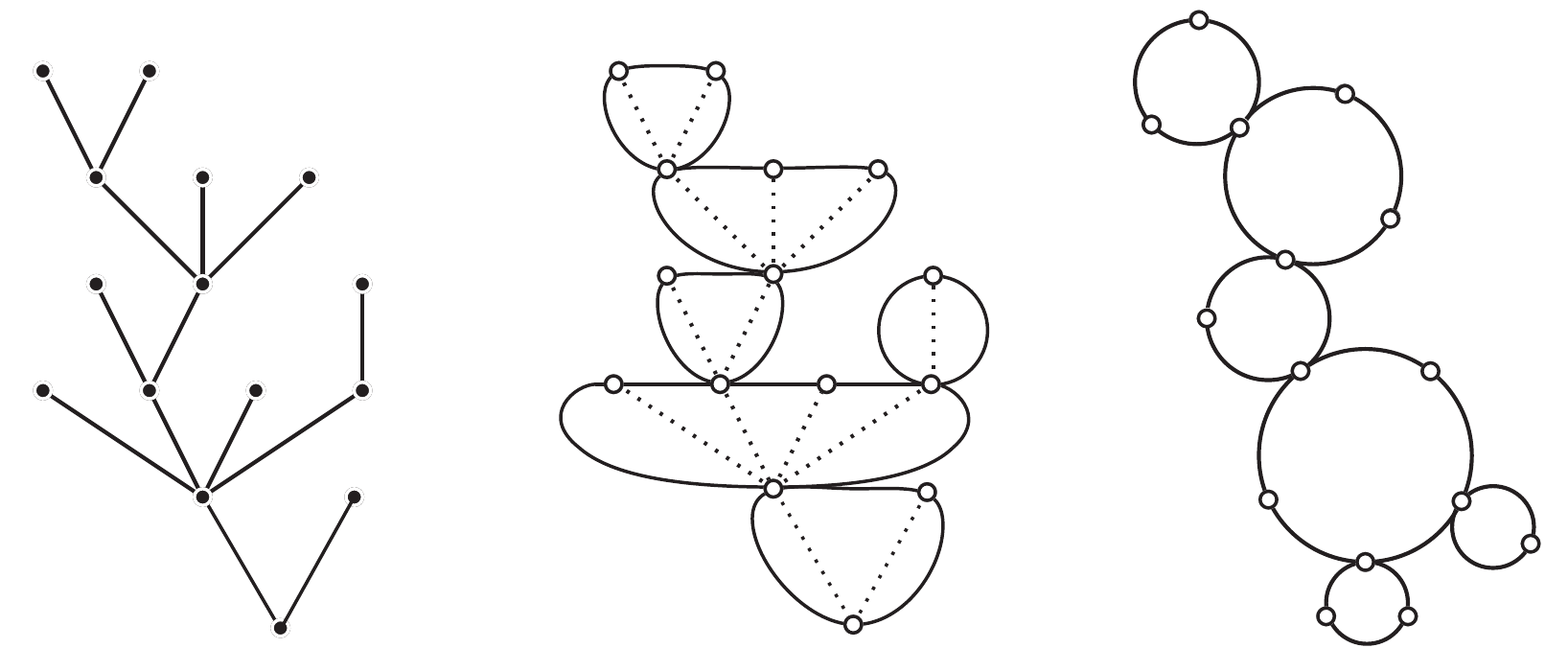}
 \caption{ \label{fig:loop}A plane tree $\tau$ (left) and its looptree $ \mathsf{Loop}( \tau)$ (middle and right).}
 \end{center}
 \end{figure}

\bigskip\noindent\textbf{Stable looptrees.} Random stable looptrees are random compact metric spaces and can, in a certain sense, be seen as the dual of the stable trees introduced and studied in \cite{duquesne02rtl,LGLJ98}. They are constructed in \cite{CK13} using stable processes with no negative jumps, but can also be defined as scaling limits of discrete objects: with every rooted oriented tree (or plane tree) $ \tau$, we  associate a graph, called the \emph{discrete looptree} of $ \tau$ and denoted by $ \mathsf{Loop}( \tau)$,  which is the graph on the set of vertices of $\tau$ such that two vertices $u$ and $v$ are joined by an edge if and only if one of the following three conditions are satisfied in $ \tau$: $u$ and $v$ are consecutive children of a same parent, or $u$ is the first child (in the lexicographical order) of $v$, or $u$ is the last child of $v$, see Figure~\ref{fig:loop}. Note that in \cite{CK13},  $\mathsf{Loop}( \tau)$ is defined as a different graph, and that here $\mathsf{Loop}( \tau)$ is the graph which is denoted by $\mathsf{Loop}'( \tau)$ in \cite{CK13}. We view $ \mathsf{Loop}( \tau)$ as a compact metric space by endowing its vertex-set with the graph metric (every edge has unit length).

  Fix $ \alpha \in (1,2)$.  Now let $\tau_{n}$ be a Galton--Watson tree conditioned on having $n$ vertices, whose offspring distribution $\mu$ is critical and satisfies $\mu_{k} \sim c \cdot k^{-1-\alpha}$ as $k \rightarrow \infty$ for a certain $c>0$.   In \cite[Section 4.2]{CK13}, it is shown that there exists a random compact metric space $\mathscr{L}_{\alpha}$, called the stable looptree of index $ \alpha$, such that 
\begin{equation} \label{eq:invprinc}  n^{-1/\alpha} \cdot \mathsf{Loop}( \tau_{n}) \tol \left(c\, |\Gamma(- \alpha)| \right)^{-1/ \alpha} \cdot \mathscr{L}_{\alpha}, \end{equation}
where the convergence holds in distribution for the Gromov--Hausdorff topology and where $c \cdot M$ stands for the metric space obtained from $M$ by multiplying all distances by $c >0$. 

 \begin{figure}[!ht]
 \begin{center}
  \includegraphics[height=5cm]{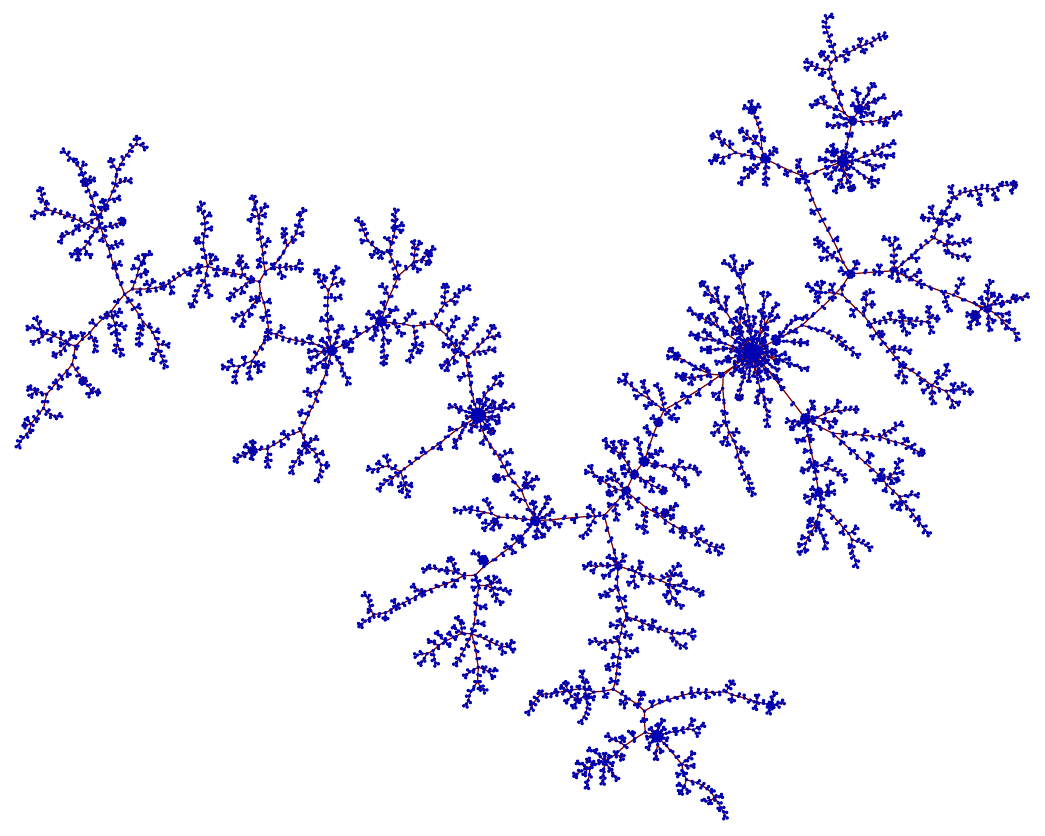}
    \includegraphics[height=4cm]{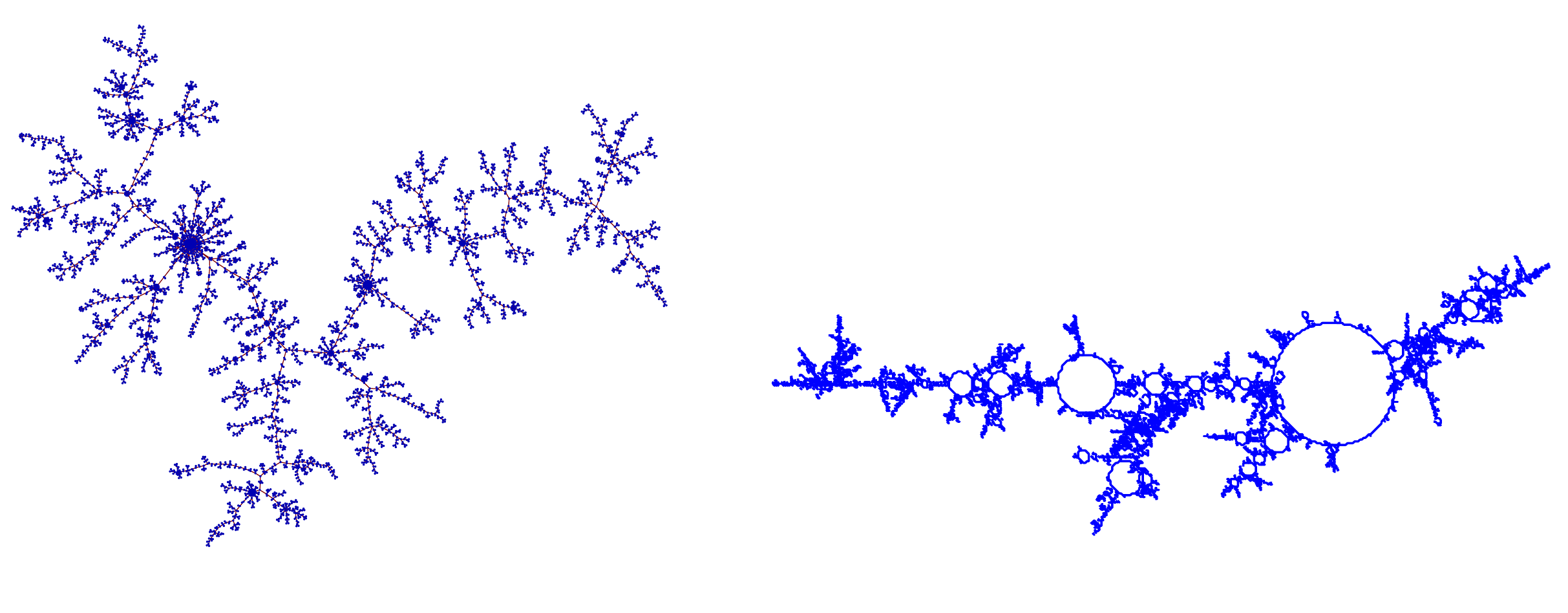}
 \caption{ \label{fig:frerots} An $ \alpha={3}/{2}$ stable tree, and its {associated}
 looptree {$\mathscr{L}_{3/2}$}, embedded in a non isometric and non proper fashion in the plane.}
 \end{center}
 \end{figure}
 
It has been proved in \cite{CK13} that the Hausdorff dimension of $ \mathscr{L}_{ \alpha}$ is almost surely equal to $\alpha$. Furthermore, the stable looptrees can be seen as random metric spaces interpolating between the unit length circle $ \mathcal{C}_{1}\de\frac{1}{2\pi} \cdot \mathbb{S}_{1}$ and  Aldous's Brownian CRT  \cite{aldous91crt} (which we view here as the real tree $ \mathcal{T}_{ \mathbf{e}}$ coded by a normalized Brownian excursion~$\mathbf{e}$, see \cite{legall05rta}). We are now in position to describe the possible scaling limits of the boundary of percolation clusters in the UIPT.  For fixed $a \in (0,1)$, let $ \partial \mathcal{H}_{a}^{\circ}(n)$ be the boundary of the white hull of the origin conditioned on the event that $ \mathcal{H}^{\circ}_{a}$ is finite and that the perimeter of $ \partial  \mathcal{H}^{\circ}_{a}$ is $n$. We view $  \partial \mathcal{H}_{a}^{\circ}(n)$ as a compact metric space by endowing its vertex-set with the graph metric.

 \begin{thm}[Curien--Kortchemski \cite{CKpercolooptrees}] \textbf{\textup{Scaling limits for $\partial \mathcal{H}^{ \circ}_{a}$ when $ \# \mathcal{H}_{a}^\circ < \infty$.}} \label{thm:scalingperco}

For every $a \in (0,1)$, there exists a positive constant $C_{a}$ such that  the following convergences hold in distribution for the Gromov--Hausdorff topology:
$$\begin{array}{cllcl} (i) &  \mbox{when }1/2<a<1,&   n^{-1} \cdot \partial \mathcal{H}_a^{\circ}(n) &\xrightarrow[n\to\infty]{(d)}&  C_{a} \cdot \mathcal{C}_{1}, \\ \ \\
(ii) & \mbox{when }a=a_{c}=1/2,&  n^{-2/3} \cdot\partial \mathcal{H}_a^{\circ}(n)  &\xrightarrow[n\to\infty]{(d)}&  3^{1/3}  \cdot \mathscr{L}_{3/2}, \\ \ \\ 
(iii) & \mbox{when }0<a<1/2, &  n^{-1/2} \cdot \partial \mathcal{H}_a^{\circ}(n) &\xrightarrow[n\to\infty]{(d)}& C_a\cdot  \mathcal{T}_{ \mathbf{e}}.  \end{array}$$
\end{thm}

Although Theorem \ref{thm:scalingperco} does not imply that $1/2$ is the critical threshold for percolation on the UIPT (as shown in \cite{angel03gpu}), it is  a compelling evidence for it.
Let us give a heuristic justification for the three limiting compact metric spaces appearing in the statement of this theorem. Imagine that we condition the cluster of the origin to  be finite and have a very large, but finite, boundary. In the supercritical regime $(i)$, as soon as  the cluster grows arms it is likely to become infinite, hence the easiest way to stay finite is to look like a loop. On the contrary, in the subcritical regime $(iii)$, having a large boundary costs a lot, so  the cluster adopts the shape which maximizes its boundary length for fixed size: the tree. In the critical case $(ii)$, these effects are balanced and a fractal object emerges: not quite a loop, nor a tree, but a looptree!

Also, the proof of Theorem \ref{thm:scalingperco} gives the expression of $C_{a}$ in terms of certain quantities involving Galton--Watson trees, which allows to obtain the  near-critical scaling behaviors of $C_{a}$ as $a \uparrow 1/2$ and $a \downarrow 1/2$.

Let us mention that the exponents appearing in the previous theorems are expected to be universal, and it is believed that the techniques of~\cite{CKpercolooptrees} may be extended to prove that the stable looptrees $(\mathscr{L}_{\alpha} : \alpha \in (1,2))$ give the scaling limits of the outer boundary of  clusters of suitable statistical mechanics models on random plane triangulations.

\section{On the distance-profile of rooted simple maps [by G. Collet]}

In~\cite{BCF14}, we prove the convergence of the distance-profile for rooted \emph{simple} maps, that is, plane maps with neither loops nor multiple edges. We also show that it implies the same type of result for the class of loopless maps, and for the class of all plane maps.
We now give a few definitions in view of stating our main result. 
For a rooted plane map $\m$ with $n$ edges, the \emph{distance} $d(e)$ 
of an edge~$e$ (with respect to the root) 
is the length of a shortest path starting at (an extremity of) $e$ and ending at the root-vertex;
the \emph{distance-profile} of $\m$ is the $n$-set $\{d(e)\}_{e\in E}$, where~$E$ denotes the edge-set of~$\m$. (Notice that we consider a distance-profile at edges, not at vertices.) 
 Let us now give some terminology 
for the type of convergence results to be obtained. We denote by $\Mone$ 
 the set of probability measures on $\mathbb{R}$, endowed with the \emph{weak topology} (that is, the topology given by the convergence in law). 
For $\mu\in\Mone$, we denote by $F_{\mu}(x)$ the cumulative function of $\mu$, 
$$\mathrm{inf}(\mu)\de\mathrm{inf}\{x:\ F_{\mu}(x)>0\}\qquad\text{ and }\qquad \mathrm{sup}(\mu)\de\mathrm{sup}\{x:\ F_{\mu}(x)<1\},$$
and we define the \emph{width} of $\mu$
as $\mathrm{sup}(\mu)-\mathrm{inf}(\mu)$. 
We also define the \emph{nonnegative shift} of $\mu$ as the probability measure (with support in $\mathbb{R}_+$) whose cumulative function is $x\mapsto F_{\mu}(x+\mathrm{inf}(\mu))$. 


\begin{defi}
We let $\muISE$ be a random variable with ISE law (as defined at the end of Section~\ref{secbs}), and $\muISEsh$ be its non-negative shift; these are random variables taking values in $\Mone$.
A sequence $\mu^{(n)}$ of random variables taking values in $\Mone$ is said to satisfy the \emph{ISE limit property} if the following properties hold:

\begin{itemize}
\item $\mu^{(n)}$ converges in law to $\muISEsh$ (for the weak topology on $\Mone$).  
\item $\sup(\mu^{(n)})$ converges in law to $\sup(\muISEsh)$ (i.e., the width of $\muISE$). 
\end{itemize}
\end{defi}

For $\mu\in\Mone$, we denote by $X(\mu)$ a real random variable with distribution given by $\mu$.
It is easy to see that if a sequence $\mu^{(n)}$ of random variables taking values in $\Mone$
converges in law to $\mu$, then $X(\mu^{(n)})$ converges in law to $X(\mu)$. It is known
that $X(\muISEsh)$ is distributed as $\mathrm{sup}(\muISE)$ (whose cumulative function has an explicit expression, see~\cite{BDFG03}). 
Hence, if  $\mu^{(n)}$ has the ISE limit property, then $X(\mu^{(n)})$
converges in law to  $\mathrm{sup}(\muISE)$. 

For an $n$-set $\mathbf{x}=\{x_1,\ldots,x_n\}$ of nonnegative values, and for $a>0$, 
define $\mu_{a}(\mathbf{x})$ as the probability measure
$$
\mu_{a}(\mathbf{x})=\frac{1}{n}\sum_{i=1}^n \delta_{x_i/(an)^{1/4}}, 
$$
where $\delta_x$ denotes the Dirac measure at $x$. Our main result is the following:

\begin{thm}[Bernardi--Collet--Fusy \cite{BCF14}]\label{theo:profile}
For $n\geq 1$, let $\pi^{\rm (s)}_n$ be the distance-profile of a  
 uniformly random rooted simple map with $n$ edges. 
 Then $\mu_{2}(\pi^{\rm (s)}_n)$ satisfies the ISE limit property. 
\end{thm}

It has been shown in~\cite{GaoWor99,BaFlScSo01} that when $n$ gets large, a uniformly random rooted loopless map of size $n$ has almost surely a unique ``giant component,'' which is a uniformly random simple map whose size is concentrated around $2n/3$, and the second largest component has size $O(n^{2/3+\delta})$ for any $\delta>0$.
Hence the following can be deduced from Theorem~\ref{theo:profile}:


\begin{corol}[Bernardi--Collet--Fusy \cite{BCF14}]
For $n\geq 1$, let $\pi^{\rm (l)}_n$ be the distance-profile of a  
 uniformly random rooted loopless map with $n$ edges. 
 Then $\mu_{4/3}(\pi^{\rm (l)}_n)$ satisfies the ISE limit property. 
\end{corol}

Similarly, it has been shown in~\cite{GaoWor99,BaFlScSo01} that when $n$ gets large, a uniformly random rooted map of size $n$ has almost surely a unique ``giant component'' $L$, which is a uniformly random loopless map whose size is concentrated around $2n/3$, and the second largest component 
has size $O(n^{2/3+\delta})$ for any $\delta>0$. 
Hence we recover here a known result, which alternatively follows from the study by~\cite{chassaing04rpl} of the profile of random rooted quadrangulations, 
combined with the recent profile-preserving bijection in \cite{ambjornbudd} between quadrangulations and maps: 

\begin{corol}[Bernardi--Collet--Fusy \cite{BCF14}]
For $n\geq 1$, let $\pi^{\rm (g)}_n$ be the distance-profile of a  
 uniformly random rooted general map with $n$ edges. 
 Then $\mu_{8/9}(\pi^{\rm (g)}_n)$ satisfies the ISE limit property. 
\end{corol}

\bibliographystyle{alpha}
\bibliography{_biblio}
\end{document}